\documentclass[12pt]{amsart}
\usepackage{amssymb}
\usepackage{amsmath} 
\usepackage{amsthm} 
\usepackage{calc} 
\usepackage{verbatim} 
\usepackage[dvips]{graphics}
\usepackage[dvips]{epsfig}
\newtheorem{thm}{Theorem}[section]
\newtheorem{lem}[thm]{Lemma}

\newtheorem{cor}[thm]{Corollary}

\newtheorem{eg}[thm]{Example}
\newtheorem{prop}[thm]{Proposition}

\newcounter{probno}
\stepcounter{probno}

\renewcommand{\qed}{\hfill$\Box$ \par \medskip}

\newcommand{\ignore}[1]{}

\newcommand{\self}{\circlearrowleft}

\newcommand{\eqdef}{:=}

\newcommand{\norm}[2][{}]{\left\|#2\right\|_{#1}}

\newcommand{\R}{\mathbf{R}}

\newcommand{\C}{\mathbf{C}}
\newcommand{\cp}{\mathbf{P}}
\newcommand{\Z}{\mathbf{Z}}
\newcommand{\N}{\mathbf{N}}

\newcommand{\supp}{\mathrm{supp}\,}

\newcommand{\pic}{\mathrm{Pic}\,}

\input epsf.sty

\renewcommand{\eqdef}{:=}
\newcommand{\as}{AS }
\newcommand{\res}{\operatorname{res}}

\title{Invariant curves for birational surface maps}

\date{\today}

\author{Jeffrey Diller \& Daniel Jackson \& Andrew Sommese}
\address{Department of Mathematics\\
         University of Notre Dame\\
         Notre Dame, IN 46556}
\email{diller.1@nd.edu, djackso1@nd.edu, sommese@nd.edu}

\thanks{The authors gratefully acknowledge support from the National
Science Foundation during the preparation of this paper.}
\begin{document}
\begin{abstract}
We classify invariant curves for birational surface maps that are expanding
on cohomology.  When the expansion is exponential, the arithmetic genus
of an invariant curve is at most one.  This implies severe constraints on both
the type and number of irreducible components of the curve.  In the 
case of an invariant curve with genus equal to one, we show that there 
is an associated invariant meromorphic two form.  
\end{abstract}
\maketitle
\markboth{\today}{\today}

\section{Introduction}
\label{intro}
One of the first things (see e.g. \cite{Bea91}, p65) one learns about dynamics 
on the Riemann sphere is that no non-trivial rational function $f:\cp^1\self$ 
leaves more than two points totally invariant.  This fact, an elementary 
consequence of the Riemann-Hurwitz Theorem, has been generalized
\cite{FoSi92}, \cite{SSU00}, \cite{BCS04} to holomorphic maps $f:\cp^k\self$ 
in any dimension.  For example, if $f:\cp^2\self$ is holomorphic, then the 
largest totally $f$-invariant curve $C = f^{-1}(C)$ is a union of at most 
three lines.  In this paper, we classify curves invariant under a 
non-trivial birational map of two complex variables.

There are two important differences between the birational and holomorphic 
cases.  First of all,
it is much too restrictive to consider only those curves which are totally
invariant in the strictest sense.  If $f:\cp^2\self$ is a birational map that
is not linear, then no curve $C$ can be $f$-invariant if one allows 
components of the critical set of $f$ in the preimage.  So instead, we define 
the preimage $f^{-1}(C)$ of a curve $C$ to exclude the critical set, and we
say that $C$ is \emph{invariant} if $C=f^{-1}(C)$.  That is, $C$ is equal
to its proper transform under $f$.

The second way in which the birational case is different is that one cannot
hope to obtain a general theorem concerning birational maps of $\cp^2$
without also considering birational maps of other complex surfaces.  The
problem is as follows.  If one begins with, say, a linear map $L:\cp^2\self$
that preserves a line $\ell$ and then conjugates with a cremona transformation
$g:\cp^2\self$ of very high algebraic degree $d(g)$, then the birational map
$f\eqdef g\circ L \circ g^{-1}:\cp^2\self$ will leave invariant the 
curve $g(\ell)$, which has degree $d(g)$.  Hence there is
no limit on the degree of a birationally invariant curve. 

We observe, however, that the map $f$ in such an example always degenerates 
when one starts to iterate it.  In particular, the algebraic degree 
$d(f^n)$ of $f^n$ is not $d(f)^n$ as one would hope.  Rather, it
is bounded above by (and in typical cases equal to) $d(g)^2$.  In a sense,
both the map $f$ and the invariant curve $g(\ell)$ are relatively simple
objects disguised as something more complicated by a poor birational choice of
coordinate.  In particular, it was shown in \cite{DiFa01} that by blowing up
points in $\cp^2$ and lifting everything to the new complex surface, one can 
always arrive at a birational map that is 
\emph{algebraically stable}.  That is, $(f^n)^* = (f^*)^n$ where $f^*$
is the induced linear action of $f$ on $\pic(X)$.  In this case, the degree 
of $f$ is replaced by the spectral radius $\lambda(f)$ of $f^*$, i.e. the
so-called \emph{(first) dynamical degree} of $f$.  Our first main result (see
section \ref{curves}) is

\begin{thm} 
\label{mainthm}
Let $f:X\self$ be an algebraically stable birational map of a complex
projective surface with $\lambda(f) > 1$.  Then the genus of any connected 
$f$-invariant curve $C$ is zero or one.
\end{thm}

By genus here, we mean what is commonly called the \emph{arithmetic}
genus of $C$.  In particular, if $C$ is irreducible, then the Riemann surface
obtained by desingularizing $C$ is either the Riemann sphere or a torus.
And if $X = \cp^2$, then the theorem amounts to saying that the degree of
an $f$-invariant curve is at most three.  In fact, we will prove (see Section
\ref{curves}) somewhat more
than the assertion in Theorem \ref{mainthm}, classifying invariant curves
for any bimeromorphic map $f:X\self$ of any compact K\"ahler surface $X$ for 
which the sequence $\norm{f^{n*}}$ is unbounded. 

Beyond the intrinsic interest of Theorem \ref{mainthm}, we note that 
invariant curves play a decisive role in many dynamically interesting examples
of birational maps.  For instance, 
Example 4 in \cite{Fav2} concerns a birational map that restricts to a 
rotation on a particular line, and this example has turned out to be 
important for testing the limits of
what is known about ergodic theory of birational maps.  In another direction,
the papers \cite{BeDi05b} and \cite{BeDi06} give detailed descriptions of the 
real dynamics of some families of birational maps, and the analysis in these
papers depends heavily on the fact that indeterminacy orbits of the maps 
are constrained to lie in invariant curves.  More generally, it is natural
to consider the class of birational maps on a surface that preserve a given 
meromorphic two form.  The support of the divisor of such a two form will 
necessarily be invariant in some sense.  We show in section \ref{elliptic} that

\begin{thm} 
\label{mainthm2}
Let $f:X\self$ be an algebraically stable birational map of a complex
projective surface, and $C$ a connected $f$-invariant curve of genus one.
By contracting curves in $X$, one can arrange additionally that $-C$ is the 
divisor of a meromorphic two form $\eta$ satisfying $f^*\eta = c\eta$.
The constant $c\in\C$ is determined solely by the curve $C$ and the induced 
automorphism $f:C\self$.
\end{thm}

\noindent 
Among other things, this theorem allows us to be quite precise about
the possibilities for a genus one invariant curve (see Corollary
\ref{ellipticlistcor}).

The plan of the paper is as follows.  Section \ref{background} contains
(mostly well-known) definitions and results concerning geometry of 
surfaces and birational maps on surfaces.  It also reviews a classification
of birational self-maps from \cite{DiFa01} that we will rely on heavily.
Section \ref{curves} presents, among other things, the proof of Theorem 
\ref{mainthm}.  The central ingredient is Corollary \ref{bignefcor} which 
says that if $C$ is an invariant curve and $K_X$ the canonical class of $X$, 
then $(C+K_X)\cdot\theta \leq 0$ for some nef class $\theta\in H^{1,1}_\R(X)$. 
Section \ref{elliptic} concerns the case of genus one invariant curves and 
contains 
the proof of Theorem \ref{mainthm2}.  Section \ref{components} discusses the
number of irreducible components of an invariant curve.  In the genus one
case, we give an upper bound that depends only on the surface and not the
map.  In the genus zero case, we give an upper bound that is more
complicated, but aside from one exceptional situation, the bound again depends 
only on the surface.  


\section{Background}
\label{background}
\subsection{Complex surfaces}
Throughout this paper, $X$ will denote a complex surface, by which we mean a
connected compact complex manifold of complex dimension two.  Usually
$X$ will be rational.  The book \cite{BHPV} is a good  
general reference for complex surfaces.  Here we recount only needed facts.

Given divisors,
$D,D'$ on $X$, we will write $D\sim D'$ to denote linear equivalence,
$D\leq D'$ if $D' = D + E$ where $E$  is an effective divisor and 
$D\lesssim D'$ if $D' - D$ is linearly equivalent to an effective divisor.
By a \emph{curve} in $X$, we will mean a reduced effective divisor.  We let 
$\pic(X)$ denote the Picard group on $X$, i.e. divisors modulo linear 
equivalence.  We let $K_X\in\pic(X)$ denote the (class of) a canonical divisor 
on $X$, which is to say, the divisor of a meromorphic two form on $X$.

Taking chern classes associates each element of $\pic(X)$ with a cohomology 
class in $H^{1,1}(X)\cap H^2(X,\Z)$.  We will have need of the larger group 
$H^{1,1}_\R(X) \eqdef H^{1,1}(X)\cap H^2(X,\R)$.  We call a class 
$\theta\in H^{1,1}_\R(X)$ \emph{nef} if $\theta^2 \geq 0$ and 
$\theta\cdot C\geq 0$ for any complex curve.  We will repeatedly rely on the
following consequence of the Hodge index theorem.
\begin{thm}
\label{hodgethm}
If $\theta\in H_\R^{1,1}(X)$ is a non-trivial nef class, and $C$ is a curve, 
then $\theta\cdot C = 0$ implies that either 
\begin{itemize}
\item the intersection form is negative when restricted to divisors supported 
  on $C$; or 
\item $\theta^2=0$ and there exists an effective divisor $D$ supported on $C$ 
      such that $D \sim t\theta$ for some $t>0$. 
\end{itemize}
In particular if $\theta$ has positive self-intersection, then the 
intersection form is negative definite on $C$.
\end{thm}

\begin{proof}
The hypotheses imply that $\theta\cdot D = 0$ for every divisor $D$ 
supported on $C$.  Suppose that the intersection form restricted to $C$
is not negative definite.  That is, 
there is a non-trivial divisor $D$ with $\supp D\subset C$
and $D^2 \geq 0$.  Then we may write $D = D_+ - D_-$ as a difference
of effective divisors supported on $C$ with no irreducible components in
common.  Since $D_+\cdot D_-\geq 0$, we have
$$
0\leq D^2 \leq D_+^2 + D_-^2,
$$
so replacing $D$ with $D_+$ or $D_-$ allows us to assume that $D$ is
effective.  In particular, $D$ represents a non-trivial class in 
$H^{1,1}_{\R}(X)$.  Since $D\cdot \theta = 0$ and $\theta^2,D^2 \geq 0$, 
we see that the intersection form is non-negative on the subspace of 
$H^{1,1}_\R(X)$ generated by $D$ and $\theta$.  By Corollary 2.15 
in \cite[page 143]{BHPV}, such a subspace must be one-dimensional.
Thus $D = t\theta$ for some $t>0$. 
\end{proof}

By the \emph{genus} $g(C)$ of a curve $C\subset X$, we will mean the quantity 
$1-\chi(\mathcal{O}_C)$, or equivalently, $1+h^0(K_C)$ minus the number of 
connected components of $C$.  If $C$ is smooth and irreducible, then $g(C)$ is 
just the usual genus of $C$ as a Riemann surface.  If $C$ is merely
irreducible, then $g(C)$ is usually called the \emph{arithmetic genus} of $C$, 
and in this case it dominates the genus of the Riemann surface obtained by
desingularizing $C$.  If $C$ is connected then $g(C) \geq 0$, but our notion
of genus is a bit non-standard in that we do not generally require 
connectedness of $C$ in what follows.  For any curve $C$, connected or not, we 
have the following \emph{genus formula} 
\begin{equation}
\label{genusformula}
g(C) = \frac{C\cdot(C+K_X)}{2} + 1.
\end{equation}

\subsection{Birational maps}
Now suppose that $Y$ is a second complex surface and $f:X\to Y$ is a 
birational map of $X$ onto $Y$.  That is, $f$ maps some Zariski open subset of 
$X$ biholomorphically onto its image in $Y$.  In general the complement of 
this subset will consist of a finite union of rational curves collapsed by 
$f$ to points, and a finite set $I(f)$ of points on which $f$ cannot be 
defined as a continuous map.  We call the contracted curves \emph{exceptional} 
and the points in $I(f)$ \emph{indeterminate} for $f$.  The birational inverse 
$f^{-1}:Y\to X$ of $f$ is obtained by inverting $f$ on the Zariski open set 
where $f$ acts biholomorphically.  Note that what we call a birational map is 
perhaps more commonly called a birational correspondence, the former term 
often being understood to mean that $I(f)=\emptyset$.

We adopt the following conventions concerning images of proper subvarieties
of $X$.  If $C\subset X$ is an irreducible curve, then $f(C)$ is defined to be
$\overline{f(C-I(f))}$, which is a point if $C$ is exceptional for $f$ and a 
curve otherwise.  If $p\in X$ is a point of indeterminacy, then $f(p)$ will 
denote the union of $f^{-1}$-exceptional curves that $f^{-1}$ maps to $p$.  We 
apply the same conventions to images under $f^{-1}$. 

Our convention for the inverse image of an irreducible curve extends by 
linearity to a \emph{proper transform} action $f^\sharp D$ of $f$ on divisors 
$D$, provided we identify points with zero.  We also have the 
\emph{total transform} action $f^* D$ of $f$ on divisors obtained by pulling 
back local defining functions for $D$ by $f$.  Total transform has the 
advantage that it preserves linear equivalence and therefore descends to a 
linear map $f^*:\pic(Y)\to \pic(X)$.  We denote the proper and 
total transform under $f^{-1}$ by $f_*$ and $f_\sharp$, respectively.

In general, $f^* D - f^\sharp D$ is an effective divisor with support equal to
a union of exceptional curves mapped by $f$ to points in $\supp D$.  It will 
be important for us to be more precise about this point.  To do so, we use the 
`graph' $\Gamma(f)$ of $f$ obtained by minimally desingularizing the variety
$$
\overline{\{(x,f(x))\in X\times Y:x\notin I(f)\}}.
$$
We let $\pi_1:\Gamma(f)\to X$, $\pi_2:\Gamma(f)\to Y$ denote projections onto 
first and second coordinates.  Thus $\Gamma(f)$ is an irreducible complex 
surface and $\pi_1,\pi_2$ are \emph{proper modifications} of their respective 
targets, each holomorphic and birational and therefore each equal to a finite 
composition of point blowups.  One sees readily that 
$f = \pi_2\circ \pi_1^{-1}$, and that the exceptional and indeterminacy sets 
of $f$ are the images under $\pi_1$ of the exceptional sets of $\pi_2$ and 
$\pi_1$, respectively.  Given a decomposition $\sigma_n\circ\dots
\circ\sigma_1$ of $\pi_2$ into point blowups, we let $E(\sigma_j)$ denote the 
center of the blowup $\sigma_j$ and 
$$
\hat E_j(f) = \sigma_1^*\dots \sigma_{j-1}^* E(\sigma_j),\quad
E_j(f) = \pi_{1*}\hat E_j(f).
$$
In particular, $\bigcup \supp E_j(f)$ is the exceptional set.  We call
the individual divisors $E_j(f)$ the \emph{exceptional components} of $f$ 
and call their sum $E(f)\eqdef \sum E_j(f)$ the \emph{exceptional divisor}
of $f$.  It should be noted that, as we have defined them, the exceptional 
components of $f$ are
connected, but in general they are neither reduced nor irreducible.   

The following proposition assembles some further information about the
exceptional components.  These can be readily deduced from well-known facts 
about point blowups.  We recall that the \emph{multiplicity} 
of a curve $C$ at a point $p$ is just the minimal multiplicity of the
intersection of $C$ with an analytic disk meeting $C$ only at $p$. 

\begin{prop} 
\label{exceptional}
Let $\sigma_j$, $E_j(f)$, and $E(f)$ be as above,
and $C\subset X$ be a curve. 
\begin{itemize}
\item $E(f) = (f^*\eta) - f^*(\eta)$ for any 
meromorphic two form $\eta$ on $X$ (here $(\eta)$ denotes the divisor of
$\eta$).  Less precisely, $E(f)\sim K_X - f^* K_X.$
\item $E_j(f)$ and $E_i(f)$ have irreducible components in common if and only 
if $f(E_j(f)) = f(E_i(f))$.  If this is the case, then $i \leq j$ implies that 
$E_i(f)\leq E_j(f)$.
\item The multiplicity with which an irreducible curve $E$ occurs in 
$E(f)$ is bounded above by a constant that depends only on the number 
of exceptional components $E_j(f)$ that include $E$.  
\item $f^* C - f^\sharp C = \sum c_j E_j(f)$ where $c_j$ is the multiplicity
of $(\sigma_n\circ \dots\circ\sigma_{j+1})^\sharp(C)$ at the point 
$\sigma_j(E(\sigma_j))$. 
\item In particular, $c_j$ vanishes if $p_j \eqdef f(E_j(f)) \notin C$,
$c_j\leq 1$ if $p_j$ is a smooth point of $C$, and $c_j>0$ if $p_j\in C$
and $E_j(f)$ is not dominated by any other exceptional component of $f$.
\item Hence (in light of the 2nd and 5th items),
$\supp f^* C - f^\sharp C = f^{-1}(C\cap I(f^{-1}))$.
\end{itemize}
\end{prop}

We will also need the following elementary fact.  

\begin{lem}
\label{sing1}
Let $C\subset X$ be a curve such that no component of $C$ is exceptional
for $f$.  If $p\in C-I(f)$, then multiplicity of $f(C)$ at $f(p)$ is
no smaller than that of $C$ at $p$.  In particular, $f(p)$ is singular
for $f(C)$ if $p$ is singular for $C$.
\end{lem}

\ignore{\begin{proof}
If $f$ acts biholomorphically at $p$, the result is obvious.  Otherwise
$f$ is holomorphic near $p$ and decomposes locally as a composition of point 
blowups.  Hence it suffices to verify the lemma for the case where $f$
is a point blowup.  In this case, the result follows from the facts that
the multiplicity of the intersection of $C$ with a generic smooth disk will
be minimal and that the image of a generic smooth disk under a point blowup
will be smooth.
\end{proof}}

\subsection{Classification of birational self-maps}
Supposing that $f:X\self$ is a birational self-map, we now recall some 
additional information from \cite{DiFa01}.  First of all, there are pullback 
and pushforward actions $f^*,f_*:H^{1,1}_\R (X)\self$ compatible with the 
total transforms $f^*,f_*:\pic(X)\self$.  The actions are adjoint with respect
to intersections, which is to say that 
\begin{equation}
\label{adjoint}
f^*\alpha \cdot \beta = \alpha\cdot f_*\beta,
\end{equation} 
for all $\alpha,\beta\in H^{1,1}_\R(X)$.  Less obviously, $f^{n*}$
is `intersection increasing', meaning 
$$
(f^{n*}\alpha)^2 \geq \alpha^2
$$
The \emph{first dynamical degree} of $f$ is the quantity
$$
\lambda(f) := \lim_{n\to\infty} \norm{f^{n*}}^{1/n} \geq 1.
$$
It is less clear than it might seem that $\lambda(f)$ is well-defined, as it 
can happen that $(f^n)^* \neq (f^*)^n$ for $n$ large enough.  However, 
$\lambda(f)$ can be shown to be invariant under birational change of
coordinate and one can take advantage of this to choose a good surface on 
which to work.

\begin{thm}
\label{asthm}
The following are equivalent for a birational map $f:X\self$ on a complex
surface.
\begin{itemize}
\item $(f^n)^* = (f^*)^n$ for all $n\in\Z$.
\item $I(f^n)\cap I(f^{-n}) = \emptyset$ for all $n\in\N$.
\item $f^n(C) \notin I(f)$ for any $f$-exceptional curve $C$.
\item $f^{-n}(C)\notin I(f^{-1})$ for any $f^{-1}$ exceptional curve $C$.
\end{itemize}
By blowing up finitely many points in $X$, one can always arrange that
these conditions are satisfied.
\end{thm}

We will call maps satisfying the equivalent conditions of this theorem 
\as (for \emph{algebraically} or \emph{analytically stable}).
If $f$ is \as, then $\lambda = \lambda(f)$ is just the spectral 
radius of $f^*$.  If $X$ is K\"ahler, then there is a nef class $\theta^+$ 
satisfying
$$
f^* \theta^+ = \lambda\theta^+.
$$
From \eqref{adjoint}
we have that $\lambda(f^{-1}) = \lambda(f)$,
so we let $\theta^-$ denote the corresponding class for $f^{-1}$.  The
following theorem summarizes many of the main results of \cite{DiFa01}, and
we will rely heavily on it here.

\begin{thm}
\label{classthm}
If $f:X\self$ is an \as birational map of a complex K\"ahler surface $X$ 
with $\lambda(f)=1$, then exactly one of the following is true (after 
contracting curves in $\supp E(f^n)$, if necessary).
\begin{itemize}
\item $\norm{f^{n*}}$ is bounded independent of $n$, and $f$ is an 
automorphism some iterate of which is isotopic to the identity.
\item $\norm{f^{n*}} \sim n$ and $f$ preserves a rational fibration.  In
this case $\theta^+ = \theta^-$ is the class of a generic fiber.
\item $\norm{f^{n*}} \sim n^2$ and $f$ is an automorphism preserving an
elliptic fibration.  Again $\theta^+=\theta^-$ is the class of a generic fiber.
\end{itemize}
If, on the other hand, $\lambda(f) > 1$, then $\theta^+\cdot \theta^- > 0$
and either $X$ is rational or $f$ is (up to contracting exceptional curves) 
an automorphism of a torus, an Enriques surface, or a K3 surface.
\end{thm}

\noindent 
We remark that the classes $\theta^\pm$ are unique up to positive multiples 
whenever $\norm{f^{n*}}$ is unbounded, and indeed under the unboundedness
assumption, we have
$$
\lim_{n\to\infty} \frac{f^{n*}\theta}{\norm{f^{n*}}} = c\theta^+
$$
for any K\"ahler class $\theta$ and some constant $c = c(\theta) > 0$.
In what follows, we will largely ignore the case in which $\norm{f^{n*}}$ is 
bounded.  After all, if some iterate of $f$ is the identity map, then every
curve in $X$ will be $f$-invariant.

To close this section, we recall a result from \cite{BeDi05a}, which we
will use in section \ref{elliptic}. 

\begin{thm}
\label{criticalthm}
If $f:X\self$ is an \as birational map of a complex K\"ahler surface $X$
with $\lambda(f)>1$, then after contracting curves in $\supp E(f^n)$, we can 
arrange additionally that
$\theta^+\cdot f(p) > 0$ for every $p\in I(f)$ and 
$\theta^-\cdot f^{-1}(p) > 0$ for every $p\in I(f^{-1})$.
\end{thm}

\section{Invariant curves}
\label{curves}
Unless otherwise noted in what follows, $f:X\self$ will always be an
\as birational map on a complex K\"ahler surface $X$.
We will call a curve $C\subset X$ \emph{invariant} for $f$ if $f(C) = C$.
If $C$ is $f$-invariant, then $f$ lifts to a biholomorphism of the 
desingularization of $C$.  In particular $f$ permutes the irreducible
components of $C$ and no such component is exceptional.  Clearly $C$ is
$f$-invariant if and only if $C$ is $f^{-1}$-invariant.

\begin{thm} 
\label{sing3thm}
If $C$ is an $f$-invariant curve, then 
$C\cap I(f)$ consists of smooth points of $C$.
\end{thm}

The proof of this result boils down, essentially, to the fact that if 
$p = f(p)$ is a fixed singular point of $C$ that also lies in an exceptional 
curve, then $f^n(p)$ should eventually be much more singular for $C$ than
$p$ is, contradicting $f$-invariance of $C$.

\proof 
In order to keep the notation simpler, we will prove the equivalent statement
that $C\cap I(f^{-1})$ consists of smooth points of $C$.  We suppose in order 
to reach a contradiction that some point $p\in I(f^{-1})$
is also a singular point of $C$.  Since $f$ is \as,  
$f^k(p)$ is well-defined for all $k\geq 0$.  By Lemma \ref{sing1} and
invariance of $C$, all points in the forward orbit of $p$ are singular for
$C$.  But the singular set of $C$ is finite, so replacing $f$ by
$f^k$ if necessary, we may assume that $p$ is fixed by $f$ and that 
$C$ is the connected curve obtained by keeping only those irreducible
components containing $p$.

Because $p\in I(f^{-1})$, we have that $f^{-1}(p)$ is a connected curve
in $\supp E(f)$, and because $f(p) = p$ we see that $p\in E\subset f^{-1}(p)$ 
for some irreducible exceptional curve $E$.  
In particular, $f^* E \geq E$ by the fifth conclusion in Proposition 
\ref{exceptional}.  More generally, for $k\geq 1$, the first item in
Proposition \ref{exceptional} and the fact that $f$ is \as imply that
\begin{equation}
\label{verycrit}
E(f^k) = (f^{k*}\eta) - f^{k*}(\eta) = 
\sum_{j=0}^{k-1} f^{j*}((f^{(k-j)*}\eta) - f^*(f^{(k-j-1)*}\eta)) = 
\sum_{j=0}^{k-1} f^{j*} E(f) \geq kE.
\end{equation}
We will complete the proof by suitably interpreting \eqref{verycrit} from 
another point of view.

Let $\Gamma$ be a minimal desingularization of the graph of $f^k$ and 
$\pi_1,\pi_2:\Gamma\to X$ be projections onto source and target.  
Because $p\in I(f^{-1})$ and $f$ is \as, it follows that $p\notin I(f^k)$.  
Hence there is a neighborhood $U\ni p$ such that $\pi_1$ maps 
$\pi_1^{-1}(U)$ biholomorphically onto $U$.  We let $C_0 = \pi_1^\sharp(C)$,
$E_0 = \pi_1^\sharp E$, $p_0 = \pi_1^{-1}(p)$ and observe
that $C_0$ is a connected curve meeting $E_0$ at $p_0$.  After decomposing 
$\pi_2 = \sigma_n\circ\dots \circ \sigma_1$ into a sequence of point 
blowups, we further define 
$$
C_j = \sigma_j\circ\dots\circ\sigma_1(C_0),\quad
E_j = \sigma_j\circ\dots\circ\sigma_1(E_0),\quad
p_j = \sigma_j\circ\dots\circ\sigma_1(p_0).
$$
In particular, $C_j$ is connected and meets $E_j$ (which will be a point 
for $j$ large enough) at $p_j$.
By $f$ invariance of $p$ and $C$, we see that $C_n = C$ and $p_n = p$.

Let $m \in\N$ be larger than $g(C)$.  By the third item in Proposition 
\ref{exceptional}, we can choose $k$ large enough in \eqref{verycrit} to 
deduce that there are at least $m$-indices $j$ for which 
$p_j = E_j = \sigma_j(E(\sigma_j))$ is the point blown up by $\sigma_j$.
It is well-known (see e.g. \cite{GrHa}, p506) that blowing up a singular 
point of a curve strictly decreases its genus.  Thus,  
$g(C_{j-1}) \leq g(C_j) - 1$ for $m$ different values of $j$.  The
cumulative effect of this is that
$$
g(C_0) \leq g(C) - m < 0,
$$
contradicting the fact that the genus of a connected curve is always 
non-negative.
\qed

\begin{cor} 
\label{basicineqcor}
Let $C\subset X$ be an $f$-invariant curve.  Then up to
linear equivalence, we have
$$
f^*(C + K_X) \lesssim C + K_X 
$$
\end{cor}

\proof
In light of Theorem \ref{sing3thm} and the fifth and first items in 
Proposition \ref{exceptional}, we have
$$
f^* C \leq C + E(f) \sim C + K_X - f^*K_X,
$$
which rearranges to give the inequality we seek.
\qed

\begin{cor}
\label{bignefcor} Suppose that $\norm{f^{n*}}$ is unbounded.  Then for any 
$f$-invariant curve $C$, we have 
$$
\theta^\pm \cdot (C + K) \leq 0
$$
In particular, when $\lambda(f)>1$, $(C+K)\cdot\theta \leq 0$ 
for some nef class $\theta$ with $\theta^2>0$.
\end{cor}

\begin{proof}
Let $\theta$ be a K\"ahler class on $X$.  Then by Corollary 
\ref{basicineqcor}, we have 
$$
0 \geq \frac{1}{\norm{f^{n*}}} \theta\cdot((C+K) - f^n_*(C+K)) 
  = \frac{\theta - f^{n*}\theta}{\norm{f^{n*}}}\cdot (C+K)
  \to c\theta^+ \cdot(C+K),
$$
which verifies the first assertion.  If $\lambda(f)>1$, then Theorem 
\ref{classthm} tells us that 
$(\theta^+ + \theta^-)^2 \geq 2\theta^+\cdot\theta^- > 0$.  Hence the second 
assertion holds for the particular class $\theta^- + \theta^+$.
\end{proof}

\begin{thm}
\label{n2classthm}
Suppose that $\norm{f^{n*}}\sim n^2$.  Then any connected $f$ invariant curve
is contained in a fiber of the elliptic fibration preserved by $f$.
\end{thm}

\begin{proof}
Let $S$ be a Riemann surface and $\pi:X\to S$ be the elliptic fibration 
preserved by $f$.  Then $\theta^+$ is the class of a fiber of $\pi$.  Since 
the self-intersection 
of a fiber is zero and generic fibers are smooth elliptic curves, the genus 
formula tells us that $0 = \theta^+(\theta^+ + K_X) = \theta^+\cdot K_X$.  
Therefore, from Corollary \ref{bignefcor} we see that any $f$ invariant curve 
$C$ satisfies $C\cdot\theta^+ \leq 0$.  Since $\theta^+$ is nef, it follows
that $C\cdot \theta^+ = 0$.
This can only happen if $C$ is contained in some fiber of $\pi$.
\end{proof}

\begin{thm}
\label{nclass}
Suppose that $\norm{f^{n*}}\sim n$ and let $\pi:X\to S$ denote the rational
fibration preserved by $f$.  Let $C$ denote the union of all irreducible 
$f$-invariant curves not contained in fibers of $\pi$.  If $C$ is non-empty,
then exactly one of the following is true.
\begin{itemize}
\item $C$ consists of one or two irreducible components, each mapped
biholomorphically by $\pi$ onto $S$.
\item $C$ consists of one irreducible component, and $\pi:C\to S$ is a two
to one branched cover.  In this case, some power of $f$ induces the identity
map on $S$.
\end{itemize} 
\end{thm}

In particular if the surface $X$ is rational, then the base of the rational
fibration is $\cp^1$.   Therefore, $f$-invariant curves must be rational or,
in the case where some iterate of the induced map is trivial, hyperelliptic.
We observe that hyperelliptic curves really do arise in this fashion.  If, 
for example, $A$ and $B$ are meromorphic functions on $\cp^1$ and 
$f:\cp^1\times\cp^1\self$ is given by
$$
f(x,y) = \left(x,-\frac{A(x)y+B(x)}{y}\right),
$$
then $f$ restricts to the identity map on the curve 
$y^2 + A(x)y + B(x) = 0$.

\begin{proof}
As in the previous proof, $\theta^+$ is the class of a generic fiber of $\pi$,
and the genus formula implies that $\theta^+\cdot K_X = -2$ (this time fibers
are rational rather than elliptic curves).  Therefore if $C$ is an
$f$-invariant curve, we see from Corollary \ref{bignefcor} that 
$C\cdot \theta^+ \leq 2$.  After removing all components of $C$ contained in 
fibers of $\pi$, we obtain that $\pi|C$ is an at most 2 to 1 branched cover of 
$S$ semiconjugating $f|C$ to an automorphism $\tilde f:\cp^1\self$.  In 
particular, $C$ has at most two irreducible components, and if there are two, 
they are necessarily isomorphic to $S$.  If $\tilde f$ is not periodic, then 
$S$ is either a torus and $\tilde f$ an ergodic translation or $S$ is $\cp^1$ 
and $\tilde f$ a linear fractional transformation.  In either case $f|C$ is 
also aperiodic.  Hence $C$ is also a torus or an elliptic curve.  In the case 
where $\tilde f$ is a linear fractional transformation, $f^n|C$ must have a 
finite, non-zero number of fixed points for all $n\in\N$ which means that $C$ 
is rational.
\end{proof}

\begin{thm}
\label{main}
Let $f:X\self$ be an \as birational map with $\lambda(f) > 1$.  Then any 
$f$-invariant curve $C$ satisfies $g(C)\leq 1$.  In particular, if $C$
is connected, then $g(C)$ is $0$ or $1$. 
\end{thm}

\begin{proof}
Let us suppose first that $X$ is irrational.  Then by Theorem \ref{classthm}, 
we
may blow down curves in $E(f^n)$ (in particular not components of $C$)
and assume that $X$ is an Enriques surface, a $K3$ surface or a torus,
and $f$ is an automorphism.  In the last two cases
$K_X = 0$, so Corollary \ref{bignefcor} tells us that $C\cdot\theta = 0$ for 
some nef class $\theta$ with positive self-intersection.  It follows that 
$C \cdot (C+K_X) = C^2 < 0$, 
so by the genus formula $C$ has genus zero and is in fact a smooth
rational curve with self-intersection -2. 
If $X$ is an Enriques surface, it is double-covered by a K3-surface
(\cite{BHPV}, page 339), 
and the theorem follows from lifting everything to the $K3$ surface.  

Now we turn to the rational case.
With $V$ as in the theorem and $K_X,K_C$ the canonical bundles of $X$ and $C$,
respectively, we have the standard short exact sequence of line bundles
$$
0 \to K_X \to K_X + C \to K_C \to 0,
$$
which gives rise to the long exact sequence
$$
\dots \to H^0(X,K_X + C) \to H^0(C,K_C) \to H^1(X,K_X) \to \dots
$$
However, since $X$ is a rational surface, $H^1(X,K_X)$ vanishes, 
so the sections of $K_X + C$ surject onto those of $K_C$.  But Corollary
\ref{bignefcor} tells us that $K_X + C$ has non-positive intersection with a 
nef class with positive self-intersection.  It follows then from Theorem
\ref{hodgethm} that 
$$
0 \leq h^0(K_C) \leq \dim h^0(K_X+C) \leq 1. 
$$ 
We conclude 
$$
g(C) = 1 - h^0(\mathcal{O}_C) + h^0(K_C) \leq h^0(K_C) \leq 1.
$$
\end{proof}


\section{Genus 1}
\label{elliptic}
For the remainder of this paper, we will assume $f:X\self$ is a birational
map on a complex K\"ahler surface, as in the hypotheses of Theorem \ref{main}.
Specifically, $f$ is \as and $\lambda(f) > 1$.  For irreducible invariant 
curves $C$, the implications of Theorem \ref{main} are clear:  $C$ is either a 
rational curve with at
most one simple cusp or normal crossing, or $C$ is a smooth elliptic curve.
But the conclusion of the theorem also says much about reducible invariant
curves.  For instance, in the genus 0 case, the following is well-known
(see e.g. \cite{BHPV}, page 85)

\begin{thm}
\label{genus0} A connected curve $C=C_1\cup\dots\cup C_k$ has genus 0 if and 
only if 
$C$ is a tree of smooth rational curves $C_j$.  That is, 
$C'\cdot C'' = 1$ for every decomposition of $C = C'\cup C''$ into 
connected curves without common components.
\end{thm}

In the case where $C$ is a connected genus 1 invariant curve, it turns out
that we can give an even more specific description.  To do this, it
is necessary first to refine the argument from Theorem \ref{main} and
prove

\begin{thm}
\label{genus1}
Let $f:X\self$ be an \as birational map with 
$\lambda(f) > 1$, and suppose that $V = f(V)$ is a connected invariant curve 
with $g(V) = 1$.  Then by contracting finitely many curves,
one may further arrange the following.
\begin{itemize}
\item $V \sim -K_X$ is an anticanonical divisor.
\item $I(f^n)\subset V$ for every $n\in\Z$.
\item Any connected curve strictly contained in $V$ has genus zero.
\item If $W$ is a connected $f$-invariant curve not completely contained in 
      $V$, then $W$ has genus zero, is disjoint from $V$, and is equal to a 
      tree of smooth rational curves, each with self-intersection $-2$.
\end{itemize}
\end{thm}

We remark that in general, contracting curves will be necessary to
achieve the conclusions of this theorem.  One can always cause the conclusions 
to fail trivially by blowing up finitely many consecutive elements in the
orbit of a point not contained in $\bigcup_{n\in\Z} I(f^n)$.

\begin{proof}  
Let us observe from the outset that contracting exceptional curves for $f$
and $f^{-1}$ will 
not disconnect $V$, nor change the fact that $f(V) = V$, and since the genus 
of $V$ cannot decrease when a curve is contracted, it follows from 
Theorem \ref{main} that $g(V)$ will remain equal to one.

After contracting exceptional curves, we can assume that the conclusion
of Theorem \ref{criticalthm} is in force.  So if $C$ is a curve containing
a point $p\in I(f^n)$, then we have that
$$
\theta^+ \cdot C = \frac{f^{n*}\theta^+}{\lambda^n}\cdot C
                 = \theta^+ \cdot \frac{f^n_* C}{\lambda^n} > 0,
$$
because the last item in Proposition \ref{exceptional} tells us that $f^{n*}C$
contains the the $f^{-n}$-exceptional curve $f^n(p)$.  Similarly, 
$\theta^-\cdot C > 0$ whenever $C\cap I(f^{-n}) \neq 0$.  

From the proof of Theorem \ref{main} and the assumption that $g(V) = 1$, we 
see that $h^0(K_X) = h^0(K_V) = 1$.  In particular, the line bundle $K_X + V$ 
is effective.  So assuming that $K_X \not\sim -V$, we have that
$$
K_X + V \sim D,
$$
where $D$ is a non-trivial effective divisor.  We will show that we can
contract a set of mutually disjoint components of $D$ without changing the 
facts that $V$ is invariant, that $g(V)=1$, and especially, that $f$ is AS.

By Corollary \ref{basicineqcor} and the
fact that the set of effective divisors is preserved by pullback
we have 
$$
0 \leq f^{n*} D \lesssim D
$$
for all $n\in\Z$.  In addition, Corollary \ref{bignefcor} implies that
every irreducible component of $D$ is intersection orthogonal to 
$\theta^+$ and $\theta^-$.  From this, we have that the intersection form is 
negative definite for divisors supported on $D$.  Thus the 
previous inequality actually holds on the level of divisors:
$$
0\leq f^{n*} D \leq D.
$$
The second paragraph of the proof implies that $I(f^n)\cap \supp D =
\emptyset$ for all $n\in\Z$.  So for each irreducible component $C$ of $D$ 
and each $n\in\Z$, we have that $f^n(C) = f^n_*C$ is either a point (i.e. 0) 
or an irreducible curve dominated by $D$.  In the latter case, the absence of
points of indeterminacy in $D$ combines with Lemma \ref{sing1} to imply that 
$f^n(C)$ will be smooth if and only if $C$ is.

The irreducible components $C$ of $D$ are now seen to fall into two classes: 
\emph{periodic components}, satisfying $C = f^n(C)$ for some $n\in\N$,
and \emph{eventually exceptional} components satisfying $C\subset\supp E(f^n)$ 
for $n\in\N$ large enough.  Since a component is periodic for $f$ if and only 
if it is periodic for $f^{-1}$, it must also be the case that a component is 
eventually exceptional for $f$ if and only if it is eventually exceptional for 
$f^{-1}$.  

To find a component of $D$ to contract, we apply the genus formula and
the hypothesis $g(V)=1$ to arrive at $D\cdot V = 0$.  Thus
$$
0 > D^2 = D\cdot (V+K_X) = D\cdot K_X.
$$  
We can therefore choose an irreducible component $C$ of $D$ satisfying 
$C\cdot K_X < 0$.  Because $C^2<0$, (\cite{BHPV}, p91) tells us that such a 
component must actually be a smooth rational curve with self-intersection 
$-1$---i.e. contractible.  Any non-trivial image $f^n(C) = f^n_* C$ of $C$ 
will therefore also be a smooth rational curve.  Because $f^n_*$ is 
intersection increasing and $(f^n_* C)^2<0$, we see that $(f^n_* C)^2$ must 
also be $-1$.  Finally, if $f^n_*C$ is distinct from $C$, then 
$(f^n_*C + C)^2 < 0$ 
implies that $f^n(C)\cap C = \emptyset$.  Hence the entire one-dimensional
portion of the orbit of $C$ can be contracted simultaneously, yielding
a smooth surface in which each irreducible curve in the orbit of $C$ has been 
replaced by a point.

The map $f$ descends to a birational map on this new surface.  In the
case where the contracted curves are eventually exceptional, there exist both
points of indeterminacy and exceptional curves that are eliminated by the 
contraction, but in neither case is a point of indeterminacy or an exceptional 
curve created.  Hence $f$ remains \as after the contraction.  As in the first 
paragraph of the proof, the connectedness, the invariance, and the genus of
$V$ are unaffected by the contraction (even though, in the case where $C$ is 
periodic, it could happen that $C\subset V$).  So after contracting, either we 
have $V\sim -K_X$, or we can repeat the preceding argument to contract yet more
curves in $X$.  The dimension of $\pic(X)$ is finite and decreases with every 
contraction, so eventually the process will end and it will then be the case 
that $V\sim -K_X$.  That is, we have established the first assertion of the 
theorem.

Since $V \sim -K_X$, the first and fifth conclusions of Proposition 
\ref{exceptional} together with Lemma \ref{sing3thm} imply for every 
$n\in\Z$ that
$$
E(f^n) + V \sim f^{n*} V \leq E(f^n) + V,
$$ 
which can only happen if $f^{n*}V = E(f^n) + V$.  The fifth item in
Proposition \ref{exceptional} therefore yields additionally
that $I(f^n) \subset V$ for all $n\in\Z$.  That is, the second assertion
in the theorem holds.  

Now suppose that $W$ is a connected curve strictly contained in $V$, and 
let $W' = V - W$ be the complementary curve.  Since $V$ is connected, 
we have $W\cdot W' > 0$.  Hence $V\sim -K_X$ implies that
$W\cdot (W+K) = W\cdot(-W') < 0$.  So by the genus formula, $g(W) = 0$,
and the third assertion is proved.

Finally, we consider a connected $f$-invariant curve $W$ not completely 
contained
in $V$.  If $W\cap V\neq \emptyset$, then $W\cup V$ is a connected 
$f$-invariant curve of genus at least one and not linearly equivalent to 
$-K_X$.  Therefore, we can apply the first assertion of the theorem to 
$V\cup W$ in place of $V$, blowing down curves in $X$ until $V\cup W$ descends 
to a curve linearly equivalent to $-K_X$.  Again, this process can only be 
repeated finitely many times, so after contracting more curves if necessary,
we can suppose with no loss of generality that an $f$-invariant curve 
$W\not\subset V$ is actually disjoint from $V$.  In particular 
$W\cdot K_X = 0$ and $W\cap I(f^n)=\emptyset$ for every $n\in\Z$.
The latter property implies that $W = f^*W = f_*W$ is invariant on the level
of total, as well as proper, transform.  Therefore, 
$$
\lambda\theta^+ \cdot W = f^*\theta^+\cdot W = \theta^+\cdot f_*W =
\theta^+\cdot W, 
$$ 
which implies that $\theta^+\cdot W = 0$.  Similarly, $\theta^-\cdot W = 0$.
Since $\theta^+ +\theta^-$ is nef and has positive self-intersection, we
conclude from Theorem \ref{hodgethm} that the intersection form restricted to 
$W$ is negative definite.  The last assertion in Theorem \ref{genus1} now
follows from the discussion of A-D-E curves in \cite[page 92]{BHPV}.  
\end{proof}

\begin{cor}
\label{ellipticlistcor}
The curve $V$ in the conclusion of the Theorem \ref{genus1} is one of the
following
\begin{itemize}
\item a smooth elliptic curve;
\item a rational curve with an ordinary cusp;
\item a rational curve with a single normal crossing;
\item a union of two smooth rational curves meeting tangentially;
\item a union of two smooth rational curves meeting transversely at
two distinct points.
\item a union of three smooth rational curves intersecting at
a single point;
\item a `cycle' $C_1,\dots, C_k$ of two or more smooth rational curves 
satisfying $C_jC_{j+1} = 1$ for all $j\geq 1$, $C_kC_1 = 1$ and $C_jC_k = 0$ 
for distinct $j,k$ otherwise.
\end{itemize}
\end{cor}

By paying a little more attention in the proof of Theorem \ref{genus1}
we can extract another piece of information that will prove useful
below.

\begin{cor}
\label{nonewsingscor}
Suppose that the curve $V$ in the conclusion of Theorem 
\ref{genus1} contains a cusp, a triple point, or a tangency.  Then
$V$ has one of these singularities even before contraction.
\end{cor}

\begin{proof}
If the corollary is false, then $V$ develops a cusp, a triple point,
or a tangency in the course of contracting a $-1$ curve $C$.  Suppose
first that $C\not\subset V$.  If $C\cdot V > 1$, then the genus
of $V$ will strictly increase when $C$ is contracted.  Since the
genus of an invariant curve cannot exceed one and attains this value
even before contraction, we see that in fact $C\cdot V \leq 1$.  So
either $C$ does not intersect $V$, or $C$ meets $V$ transversely
at some smooth point of $V$.  In neither case will contracting $V$
add to or change the singular points of $V$.  That is, no cusps,
triple points, or tangencies can be created when $C\not\subset V$.

Now suppose that $C\subset V$.  Then contracting $C$ leads to a cusp
only if $C$ is tangent to some other component of $V$ and contracting
$C$ leads to a tangency only if $C$ meets two other components of $V$
at the same point.  Hence, the only real concern is that contracting
$C$ might create a triple point.  If that happens, then we see that
$C$ meets at least three other components of $V$.  Thus
$$
C\cdot V = C\cdot C + C\cdot(V-C) \geq -1 + 3 = 2.
$$
Once again, combining this estimate with the genus formula implies that
the genus of $V$ will strictly increase when $C$ is contracted.  Since
this cannot happen, the proof is complete.
\end{proof}

If $V\sim -K_X$ then there is a meromorphic two form $\eta$ on $X$,
unique up to constant multiple, with
simple poles along $V$ and no zeroes or poles elsewhere. If $V$ is 
also $f$-invariant, then 
$$
(f^*\eta) = f^*(\eta) - E(f) = - f^* V - E(f) = - V + E(f) - E(f) = 
-V = (\eta).
$$
Hence 
\begin{equation}
\label{t1}
f^*\eta = t \eta
\end{equation}
for some $t\in\C$.  The constant $t$ can be computed by looking at $f|V$,
as we will now explain.  We recall that the explicit version of the linear
transformation $H^0(X,K_X+V)\to H^0(V,K_V)$ is the Poincar\'e
residue map (\cite{BHPV}, p66), which prescribes to each  meromorphic two 
form on $X$
with simple poles along $V$ a meromorphic one form on $V$.  If $V$
has local defining function $h$ on some open set $U\subset X$, then
we can write
$$
\eta = \frac{dh}{h}\wedge \eta' + \tau,
$$
on $U$, where $\eta'$ is a holomorphic $1$-form and $\tau$ is a holomorphic
$2$-form.  It turns out that $\res(\eta) \eqdef \eta'|V$ is independent of the 
choice of defining function.  Since $h^0(K_V) = 1$ and $f|V$ is an 
automorphism, it follows that
\begin{equation}
\label{t2}
f^* \res(\eta) = t\res(\eta)
\end{equation}
for some $t\in\C$.

\begin{prop}
The constant $t$ is the same in \eqref{t1} and \eqref{t2}.
\end{prop}  

\begin{proof}
It's enough to observe that $f^*$ (applied to meromorphic one and two forms)
commutes with the Poincar\'e residue map. This is clear from the 
definition of the Poincar\'e residue, since near any
point $p\in V-I(f^{-1})-\supp E(f^{-1})$, we have that $h\circ f$ remains a 
local defining function for $V$ near $f^{-1}(p)$ and 
$$
f^*\eta = \frac{d(h\circ f)}{h\circ f} \wedge f^*\eta' + f^*\tau,
$$
on $f^{-1}(U)$.  This shows that $\res(f^*\eta) = f^*\res(\eta)$ on a dense
open subset of $V$, and by analytic continuation, the equality holds
everywhere on $V$.
\end{proof}

\begin{cor} 
\label{invformcor}
Let $f:X\self$ be an \as birational map
preserving a curve $V$ with genus one.  Then there is a meromorphic two
form $\eta$ on $X$ satisfying $f^*\eta = t\eta$ for some $t\in\C$.  
Moreover, $t$ is a root of unity unless one of the following holds: 
\begin{itemize}
\item some rational component of $V$ has an ordinary cusp;
\item two distinct rational components of $V$ meet tangentially;
\item three distinct rational components of $V$ meet transversely at a 
single point;
\end{itemize}
\end{cor}

\begin{proof}
Let $V$ be an $f$-invariant curve with genus one.  Let us first assume that
$V$ is as in the conclusion of Theorem \ref{genus1}.  Then as described 
above, we may choose $\eta$ to be a meromorphic two form with divisor
$(\eta) = -V$, and conclude that $f^*\eta = t\eta$.  In general, it will
be necessary to contract curves in $X$ to reach this situation.  However,
if $\pi:X\to \tilde X$ is the map that accomplishes the contraction, and
$\tilde \eta$ is the invariant (up to scale) two form on $\tilde X$, then
it follows that, $\eta\eqdef\pi^*\tilde\eta$ is invariant on $X$ simply because
invariance is a pointwise condition for a two form and need only 
be verified on an open set that avoids the curves contracted by $\pi$.

Now for the conclusion concerning the specific value of $t$, Corollary
\ref{nonewsingscor} allows to assume that $V\sim -K_X$ is one of the curves 
in the conclusion of Corollary \ref{ellipticlistcor}.  We recall some
facts about meromorphic one forms in the range of the Poincar\'e residue
operator associated to a curve $V\sim -K_X$.  It is evident from the 
definition above that such a form will be holomorphic and non-zero at smooth 
points of $v$.  Less obviously, a non-trivial form in the range of $\res$
will have simple poles at any normal
crossing of $V$ and double poles at other singularities.  Let us consider
the implications of these facts in two contrasting particular cases.

If $V = C'\cup C''$ is a union of two smooth rational curves intersecting 
transversely at two distinct points, then we can choose a uniformizing 
parameter $z$
for $C'$ so that the normal crossings correspond to $z=0$
and $z=\infty$.  Thus forms in the range of the residue operator will
be multiples of $\frac{dz}{z}$ in this coordinate.  On the other hand,
the restriction of $f$ to $C'$ will either preserve or switch $z=0$ and
$z=\infty$; that is $f|C':z\mapsto az$ or $f|C':z\mapsto a/z$.  From
this, it is clear that $(f|C')^* \frac{dz}{z} = \pm \frac{dz}{z}$.  In
particular, $t = \pm 1$.

Suppose instead that $V = C'\cup C''$ is a union of two smooth rational curves
intersecting tangentially at a single point.  In this case we choose a 
parameter $z$ on $C'$ so that the intersection occurs at $z=\infty$.
Thus forms in the range of $\res$ are multiples of $dz$ and $f|C':z\mapsto az +
b$.  It follows that $t = a$.  The remaining cases in Corollary 
\ref{invformcor} can 
all be analyzed in a similar fashion with the result that $t$ is a root
of unity unless $V$ is a rational curve with a cusp, a union of two lines
meeting tangentially, or a union of three lines meeting at a single point.
\end{proof}

\begin{eg}
For $a,b,c\in\C$, the map $f=f_{abc}:\cp^2\self$ given in homogeneous 
coordinates by
$$
f_{abc}:[x,y,z]\mapsto [x(x+ay+z/b),y(x/a + y + cz),z(bx+y/c+z)]
$$
is birational and preserves the curve $\{xyz=0\}$, which is a union of 
three lines and has genus 1.  The inverse of $f_{abc}$ is the map 
$f_{a^{-1}b^{-1}c^{-1}}$.  For generic values of $a,b,c$, one can check 
using the second item in Theorem \ref{asthm} that $f$ is \as.  In this
case $\lambda(f) = 2$ is just the degree of the homogeneous polynomials
defining $f$.  The meromorphic two form, given in affine coordinates by 
$dx\wedge dy/xy$ is invariant under $f$.
\end{eg}

In his thesis \cite{Jac05}, Jackson gives examples of birational maps
$f:\cp^2\self$ with genus 1 invariant curves of each of the last
four types presented in Corollary \ref{ellipticlistcor}.  Such examples
seem plentiful.  However, we are not presently aware of an 
example of an \as  birational map $f:X\self$ with $\lambda(f) > 1$ that 
preserves an \emph{irreducible} genus 1 curve.


\section{The number of irreducible components of an invariant curve}
\label{components}
Throughout this section, we take $X$ to be a rational surface and set
$$
h^{1,1} \eqdef \dim H^{1,1}(X) = \dim \pic(X).
$$
We take $f:X\self$ to be an \as birational map with first dynamical degree 
$\lambda(f)>1$, and we suppose that $C\subset X$ is a (not necessarily
connected) $f$-invariant curve.  
It would be nice to have an upper bound for the number of components of $C$
that depended only on $X$ and not on $f$.  Our results so far come 
close to giving such a bound, but in the case where each connected component 
of $C$ has genus 0, we are not presently able to rule out the possibility
that $C$ contains arbitrarily many fibers in a rational fibration.  Consider,
for instance, the following example that shows it is possible to have at
least three such fibers invariant.

\begin{eg}
Let $f:\cp^1\times\cp^1\self$ be given in affine coordinates by
$$
f(x,y) = \left(x\frac{xy+3y-4}{4xy-6x+2},\frac{xy+3y-4}{xy^2+5y-6}\right).
$$
Then $f$ is birational with $I(f)=\{(\infty,0),(0,\infty),(1,1),(-1,2)\}$.
For any choice of three points $p_1,p_2,p_3\in I(f)$, there is a unique 
hyperbola of the form $(x-a)(y-b) = c$ passing through $p_1,p_2,p_3$,
and the four hyperbolas obtained this way constitute the exceptional set of 
$f$.  The lines $\{x=0\}$, $\{x=1\}$ and $\{x=\infty\}$ are all $f$-invariant.
Moreover, if we take the classes of a horizontal and a vertical line as
a basis for $\pic(\cp^1\times\cp^1)$, then $f^*$ is given by the matrix
$$
\left(
\begin{matrix}
2 & 1 \\ 1 & 2
\end{matrix}
\right).
$$
The map $f$ is not by itself \as:  for instance, 
$f(\{xy = -2\}) = (1,1) \in I(f^{-1})\cap I(f)$.
However, if we set $g(x,y) = L\circ f$, where $L(x,y) = (x,\frac{ay+b}{cy+d})$,
then $g$ will be \as for generic choices of $a,b,c,d$, and all of the other
properties of $f$ that we have just described will be retained for $g$.
In particular, $\lambda(g)$ will be the largest eigenvalue $3$ of $f^*$.
\end{eg}

Our aim in this section is to show that, once large unions of fibers in
a rational fibration are excluded, we have the desired upper bound on 
the number of components of $C$.  We rely on the following variation on 
a result \cite[pp. 111]{BHPV} of Zariski. 

\begin{prop} 
\label{picprop}
Let $V\subset X$ be a curve and $H\subset\pic(X)$ be the subspace generated
by divisors supported on $V$.  Suppose that the intersection form is
non-positive on $H$.  Then there is a unique choice of $k$ effective 
divisors $D_1,\dots,D_k$ with the following properties.
\begin{itemize}
\item For each $1\leq j\leq k$, $\supp D_j$ is a distinct connected component 
      of $V$. 
\item A divisor $D$ supported on $V$ has self-intersection $0$ if and
      only if $D = \sum c_j D_j$ for some $c_j\in\Z$.
\item Every linear equivalence among divisors supported on $V$ has
      the form $\sum_{j=1}^k c_jD_j \sim 0$ for some $c_j\in\Z$.  
\item If $k\leq 1$, then $V$ has at most $h^{1,1} - 1$ irreducible
      components.  
\item If $k\geq 2$, then there is a fibration $\pi:X\to \cp^1$ such that each 
      $D_j$ is a (possibly non-generic) fiber of $\pi$.  Moreover,
      $V$ has at most $h^{1,1}+k-2$ irreducible components.
\end{itemize}
\end{prop}

\begin{proof}
Suppose first that $V$ is connected and that there is a non-trivial divisor
$D$ with $\supp D\subset V$ and $D^2 = 0$.  Suppose that $E\subset V$ is an 
irreducible component of $V - \supp D$ such that 
$E\cap \supp D\neq \emptyset$.  Replacing $D$ by $-D$ if necessary, we may 
assume that $D\cdot E>0$.  Thus
$$
(D + nE)^2 = 2nE\cdot D + E^2 > 0
$$
for $n$ large enough, contradicting our assumption that the intersection
form is non-positive for divisors supported on $V$.  It follows that 
$\supp D = V$.  As in the proof of Theorem \ref{hodgethm}, we may assume
that $D$ is effective.

If $D'$ is some other non-trivial effective divisor supported on $V$ with 
$D'^2 = 0$, then we may choose $a,b\in\N$ so that $aD-bD'$ vanishes on
some irreducible component of $V$.  However, our hypothesis on $V$ and
Theorem \ref{hodgethm} imply that
$$
0 \geq (aD'-bD)^2 = 2abD'\cdot D \geq 0.
$$  
That is, $(aD - bD')^2 = 0$.  If $aD-bD'$ were non-trivial, then the
argument in the first paragraph would force $\supp aD-bD'=V$.  By
construction, this is not the case, so we conclude $aD=bD'$.
Furthermore by minimizing $b/a$, we arrive at an effective divisor $D$ 
supported on $V$ such that any other divisor $D'$ supported on $V$ with 
$D'^2 =0$ is an integer multiple of $D$.

Now let us drop the assumption that $V$ is connected.  In this case,
our arguments so far yield divisors $D_1,\dots,D_k$ satisfying the first two 
conclusions of the proposition.  Any linear equivalence involving only
divisors supported on $V$ may be written
\begin{equation}
\label{longeqn}
D_0^+ + D_1^+ +\dots + D_k^+ = 
D_0^- + D_1^- + \dots +D_k^- 
\end{equation}
where $D_j^\pm$ are effective divisors with no irreducible components in
common and satisfying
\begin{itemize}
\item $\supp D_j^\pm \subset \supp D_j$ for $j=1,\dots,k$; and 
\item $\supp D_{0,+}\cap \supp D_j = \emptyset$,
\end{itemize}
for $j=1,\dots,k$.  If $D_0^+\neq 0$, then intersecting both sides of 
\eqref{longeqn} with $D_0^+$ and applying our hypothesis on $V$ implies that 
$$
0 > D_0^{+2} = D_0^+\cdot D_0^- \geq 0.
$$
Hence $D_0^+ = 0$, and similarly $D_0^-=0$.  The same argument shows for
$j\geq 1$, that $D_j^\pm$ is a multiple of $D_j$.  Since the divisors $D_j$ 
are minimal among effective divisors with self-intersection 0, the multiples 
are integers.  This establishes the third conclusion of the proposition.

To obtain the desired upper bounds on the number of irreducible components of
$V$, we observe that $H\neq \pic(X)$ because $X$ is rational and must contain 
curves with positive self-intersection.  Hence $\dim H \leq h^{1,1}-1$.  The 
third conclusion of the proposition implies that there are no more than 
$\max\{0,k-1\}$ independent equivalences among the irreducible components of 
$V$.  Hence $V$ has at most $h^{1,1}-1+\max\{0,k-1\}$ irreducible components.

Finally, if $k\geq 2$, then Theorem \ref{hodgethm} tells us that 
$aD_1 \sim bD_2$ for some $a,b\in\N$.  Since 
$\supp D_1\cap \supp D_2 = \emptyset$, we can choose a surjective
holomorphic function $\pi:X\to \cp^1$ with divisor $(\pi) = aD_1-bD_2$.  
By Stein Factorization \cite[pp. 32]{BHPV}, we can assume that $\pi$
has connected fibers.  In particular, each of the other divisors $D_j$
must have support equal to a fiber of $\pi$, and since we have
chosen $D_j$ to be minimal, we conclude for each $j$ that some integer 
multiple of $D_j$ is linearly equivalent to the generic fiber of $\pi$.
\end{proof}

\begin{thm}
\label{treecompsthm}
Suppose that every connected component of $C$ has genus 0.  Then either
\begin{itemize}
\item $C$ has at most $h^{1,1}+1$ connected components; or
\item there is a holomorphic map $\pi:X\to\cp^1$, unique up to automorphisms
      of $\cp^1$, such that $C$ 
      contains exactly $k\geq 2$ distinct fibers of $\pi$, and $C$ has at most
      $h^{1,1}+k-1$ irreducible components.
\end{itemize}
\end{thm}

Actually, we know of no examples of a birationally invariant curve $C$ equal 
to a disjoint union of connected genus 0 curves and comprising more than 
$h^{1,1}+1$ irreducible components.

\begin{proof} 
If $C$ supports no divisors with positive self-intersection, then the 
corollary follows immediately from  Proposition \ref{picprop}.  So we
assume with no loss of generality that $C$ supports a divisor $D$ with 
positive self-intersection.  We may assume that $\supp D$ is 
connected, since $(D_1 + D_2)^2 = D_1^2 + D_2^2$ for divisors with disjoint
supports.  Hence there is a connected component $C_0\subset C$ containing
$\supp D$, and by Theorem \ref{hodgethm}, the intersection form is negative
definite on $C-C_0$.

Let $E_0\subset \supp D$ be any irreducible component.  If $C_0 - E_0$
supports no divisors with positive self-intersection, then once again the
theorem follows from Proposition \ref{picprop}.  Otherwise, we apply the
argument from the first paragraph and obtain a unique connected component
$C_1$ of $C_0 - E_0$ that supports a divisor with positive self-intersection.
Moreover, by Theorem \ref{genus0}, the curve $C_0$ is a tree
of smooth rational curves, and there is a unique irreducible component $E_1$
of $C_1$ that meets $E_0$.  We then (uniquely) continue this process, 
obtaining two finite sequences of curves $E_j,C_j\subset C$, $0\leq j\leq n$, 
subject to the following conditions.
\begin{itemize}
\item $C_j$ supports a divisor with positive self-intersection; 
\item $C_{j+1}$ is a connected subtree of $C_j$.
\item $E_j$ is the unique irreducible component of $C_j$ that meets $C-C_j$.
\end{itemize}
The terminal index $n$ is the first for which the intersection form
$C_n-E_n$ supports no divisors with positive self-intersection. 
That is, the intersection form is non-positve on $C_n-E_n$.  

There are now two possibilities to consider.  The first is that the 
intersection form is negative definite for divisors supported on $C_n-E_n$.  
Because $C_n$ supports a divisor with positive self-intersection and
$C_n \cap (C-C_n - E_{n-1}) = \emptyset$, the intersection form is also 
negative definite for divisors supported on $C - C_n - E_{n-1}$.  All
told, the intersection form is negative for divisors supported anywhere
on $C-E_n-E_{n-1}$, and the theorem follows from Proposition \ref{picprop}.

The other possibility is that $C_n - E_n$ supports divisors with
zero self-intersection.  In this case, the intersection form will
be non-positive on $C_n - E_n$.  Let $H,H'\subset\pic(X)$ be the
subspaces spanned by divisors supported on $C-E_n$ and $C$, respectively.
Then $\dim H < \dim H'$, because $H'$ contains elements with positive
self-intersection, whereas $H$ does not.  Hence the theorem follows again
from the bound on $\dim H$ given in Proposition \ref{picprop}.
\end{proof}

\begin{thm}
If the invariant curve $C$ contains a curve of genus $1$, then 
the number of irreducible components of $C$ is no larger than $h_{1,1}+2$.
\end{thm}

\begin{proof}
By Theorem \ref{mainthm} no connected component of $C$ has genus larger
than 1; and (with our notion of genus) the genus of a sum of two disjoint 
curves is less than the sum of their individual genera.  It follows that if 
$C$ contains a curve of genus 1, then that curve is connected.  

Contracting a $-1$ curve reduces $h_{1,1}$ by one, while it reduces the number 
of connected components of $C$ by at most one.  Thus we may apply Theorem 
\ref{genus1} to conclude that $C = C' + C''$ is a disjoint union of a 
connected genus one curve $C'$ as in the conclusion of Corollary 
\ref{ellipticlistcor}, and a curve $C''$ that supports only divisors with
negative self-intersection.  

Going through the list of possibilities in Corollary \ref{ellipticlistcor},
we see that removing an irreducible component $E_0$ of $C'$ leaves us with
a (possibly empty) chain $C_0$ of smooth rational curves.  This allows us
to employ the argument from Theorem \ref{treecompsthm} and find an
irreducible component $E_1\subset C_0$ such that the intersection form
is non-positive on $C_0-E_1$ and therefore, more generally, on $C-E_0-E_1$.
Moreover, the only connected components of $C-E_0-E_1$ that can support
non-trivial divisors with zero self-intersection are the (at most) two 
components of $C_0-E_1$.  By the first and second conclusions of Proposition
\ref{picprop}, it follows that there are at most two independent divisors with
zero self-intersection.  Hence $C-E_0-E_1$ contains at most $h^{1,1}$ 
irreducible components.
\end{proof}


\bibliographystyle{texfiles/mjo}
\bibliography{texfiles/refs}
\end{document}